\documentclass[12pt]{article}

\usepackage{amsmath,amssymb,amsfonts}

\usepackage{graphicx}

\newcommand{\bbR}{\mathbb{R}}
\newcommand{\p}{\mathsf{P}}
\newcommand{\e}{\mathsf{E}}
\newcommand{\wh}{\widehat}

\begin{document}

\begin{center}
{\bf \Large Statistical methods of SNP data analysis\\
\vskip0,2cm with applications}
\end{center}
\vskip0,3cm
\begin{center}
{\bf A.V.Bulinski\footnote{Dept. of Mathematics and Mechanics,
Lomonosov Moscow State University, Moscow, 119991,
Russia.}$^{,}$\footnote{LPMA, Universit\'e Pierre et Marie
Curie (Paris-6), 4, place Jussieu 75252, Paris Cedex 05, France.},
O.A.Butkovsky$^{1}$, A.P.Shashkin$^1$,
P.A.Yaskov$^{1,}$\footnote{Steklov Mathematical Institute, Gubkina str. 8, Moscow, 119991, Russia.}}
\end{center}

\vskip0,5cm \begin{abstract} Various statistical methods  important for genetic analysis
 are considered and deve\-loped. Namely, we concentrate on
the multifactor dimensionality reduction, logic regression, random
forests and stochastic gradient boosting. These methods and their
new modifications, e.g., the MDR method with "independent rule", are
used to study the risk of complex diseases such as cardiovascular
ones. The roles of  certain combinations of single nucleotide
polymorphisms and  external risk factors are examined. To perform
the data analysis concerning the ischemic heart disease and
myocardial infarction the supercomputer SKIF "Chebyshev" of the
Lomonosov Moscow State University was employed.

\vskip0,4cm


\vskip0,4cm {\it AMS 2010 classification}: 92B15, 92D10, 65C20.

\end{abstract}

\vskip0,5cm
\section{Introduction}

The detection of genetic susceptibility to complex diseases (such as
cardiovascular, oncological ones etc.) has recently drawn  much
attention  in many leading research centers, see, e.g., \cite{Bock}
and \cite{MachLearn09}. According to the forecast of the World
Health Organization (www.who.int), in 2030 the deaths related to
cardiovascular diseases will exceed 23 millions (this year about 17
millions), the oncological diseases will take the lives of more than
11 millions of our planet inhabitants and at least 2 millions of
people will be the victims of the diabetes. Thus this research
domain is important since one expects to provide for each person the
prophylactic measures and medical treatment taking into account
his/her genetic peculiarities  which increase the risk of some
diseases and protect from the  others, see, e.g., \cite{Lengauer07}.
Individual's DNA variations are typically described  in terms of
{\it single nucleotide polymorphisms $(${\rm SNP}$)$},  i.e. the
fragments of genetic code where a nucleotide change is possible. For
more details see, e.g., \cite{SNPgeneral}. The first examples of
genetically based diseases (e.g., sicklemia) were related with a
single mutation. Contrariwise many hard diseases such as diabetes,
Alzheimer's disease and others have a complex character as they can
be provoked by mutations in different parts of the DNA code which
are responsible for the formation of certain types of proteins.
Quite a number of recent studies (see, e.g., \cite{Dai},
\cite{Schwender3} and \cite{Schwender2}) support the paradigm that
the increasing risks of complex diseases can be explained by
combinations of certain SNP  whereas  separate mutations have no
dangerous effects.

Thereupon it should be mentioned that there existed
 a longstanding demand for
statistical analysis of biological and medical data. However, only
in the first part of the 20th century, due to the classical
contributions by K.Pirson, R.Fisher, H.Cram\'er, A.Kolmogorov,
N.Smirnov, A.Wald and other prominent statisticians, the essential
progress was achieved both in theory and applications. The methods
developed were sufficient, e.g., for investigation of the efficiency
of new medicaments. The situation has changed radically at the
beginning of the 2000's when the laboratory methods of DNA analysis
provided the data related to the personal human code structure. The
achievements in decoding of the human genome have led to formation
of vast data bases in the frameworks of the International Research
projects, see, e.g.,   GAW16 \cite{GAW} and HapMap \cite{HapMap}.
Note also that software engineering plays an important role in such
studies, see, e.g., \cite{Karchin} and \cite{ZhangBuc}. The cost of
genomic analysis has fallen considerably in the last 10 years,
allowing to collect large volumes of genetic data for genetic
mapping of complex diseases. However statistical problems arising
here require new methods of inference rather than classical ones.
Indeed, the modern statistical models involve huge number of
variables, parameters, hypotheses etc., while the sample size is
usually moderate  (several hundred or sometimes several thousand of
observations, see, e.g., \cite{Tabara}). The sample design   is
limited both by costs of analysis which are still high and by
difficulties due to the sample selection. In particular, the ethnic
homogeneity should be taken into account, as well as the influence
of external risk factors such as obesity, smoking etc.

To perform
reliable  statistical inference, it is necessary to apply new
powerful tools developed in high-dimensional statistics, artificial
intelligence, information retrieval, econometrics etc. Some of them
have been adapted and further generalized in numerous papers by
biostatisticians. Among the most important SNP analysis methods
are the multifactor dimensionality reduction (MDR), logic regression
(LR),  random forests (RF) and stochastic gradient boosting (SGB).
All approaches based on these methods do not impose any strong
restrictions on the dependence structure of variables under
consideration  (apart from independence and identical distribution
of observations within certain groups). Thus a broad class of
statistical models is defined and the model providing the best
out-of-sample fit is selected.

If one deals with too many parameters,  overfitting is likely to
happen, i.e. the estimated parameters depend too much on the given
sample. As a result the constructed estimators give  poor prediction
on new data. On the other hand, application of a very sophisticated
model may not capture the studied dependence structure of various
factors efficiently. However the trade-off between the model's
complexity and its predictive power allows to perform reliable
statistical inference via new model validation techniques (see,
e.g., \cite{ArlotCelisse} and \cite{Massart03}). The main tool of
model selection is the {\it cross-validation}, see, e.g.,
\cite{CrossValid}. Its idea is to estimate parameters by involving
only
a part of the sample {\it $($training
sample$)$} and afterwards use the remaining observations {\it
$($test sample$)$} to test the predictive power of the obtained estimates.
Then an average over several realizations of randomly
chosen training and test samples is taken, see \cite{Hastie10}.

There are two closely connected research directions in  genomic
statistics. The first one
is aimed
at the disease risk estimation when the
genetic portrait of a person is known (in turn this problem involves
 estimation of disease probability and classification of
genetic data into high and low risk domains, see, e.g.,
\cite{Massart03}). The second
trend
is to identify  relevant  combinations
of SNPs having the most significant pathogenic  (or, in other way,
protective) influence. Both
 directions are presented in this paper. Moreover, the authors propose further development of various
statistical methods and apply them to study of the risks of
cardiovascular diseases. For this purpose the new software
concerning the employment of the mentioned statistical methods was
designed and used.

Due to high-dimensionality of data many numerical procedures based on the above mentioned
statistical methods are very time consuming.
The authors are grateful to the Chancellor of the Lomonosov Moscow
State University (MSU) Professor V.A.Sadovnichy and to the Deputy
Director of the MSU Research Computing Center Professor
V.V.Voevo\-din for the opportunity to use the supercomputer SKIF MSU
``Chebyshev".

This investigation was started in the framework of the project
headed by Professor V.A.Tkachuk, the Dean of the Faculty of the
Fundamental Medicine of the MSU. An overview of preliminary results
of the work was presented at the International conference
``Postgenomic methods of analysis in biology, and laboratory and
clinical medicine" in the talk by Professor A.V.Bulinski (see
\cite{Bulinskiconference} and \cite{ourreport}).

\section{Methods}

We start with some definitions. Let $N$ be the the number of
patients in the  sample   and the vector
$X^j=\left(X^j_1,\dots,X^j_n\right)$ consist of   genetic   (SNP)
and external risk factors of $j$-th individual, $j=1,..., N$. Here
$n$ is the total number of  factors, and $X^j_i$ is the value of
$i$-th variable  (characterizing SNP or external factor) of $j$-th
individual. These variables are also called   {\it explanatory
variables} or {\it predictors}. If   $X_i$ stands for an SNP,   we
set
\begin{equation}
X_i=\begin{cases} 0,&\mbox{ no mutation in }i\mbox{-th SNP,
}\\
1,&\mbox{ heterozygous mutation},
\\
2,&\mbox{ homozygous mutation.} \\
\end{cases}
\end{equation}

We assume that the external risk factors also take no more than
three values, denoted by 0, 1 and 2. For example, we can specify a
presence or an absence of  obesity (or hypercholesterolemia etc.) by
values 1 and 0 respectively. If the external factor takes more
values (e.g., blood pressure), we can divide individuals into three
groups according to its values.

Further on   $X_1^j,\dots,X_m^j$ stand for genetic data and
$X_{m+1}^j,\dots,X_n^j$ for external risk factors. Let a binary
variable $Y^j$ ({\it{response variable}}) be equal to 1 for a {\it case}, i.e. whenever
$j$-th individual is diseased, and to $-1$ otherwise (that is for a {\it control}). Set
\begin{equation}\label{xi}
\xi=(\xi^1,\dots,\xi^N)\;\;\mbox{where}\;\;
\xi^j=(X^j,Y^j),\;\;\,j=1,\dots,N.\end{equation} Suppose
$\xi^1,\dots,\xi^N$ are i.i.d. discrete random vectors having the
same law as a vector $(X,Y)$ and independent of this vector. Assume
that $X=(X_1,\ldots,X_n).$ All random vectors (and random variables)
are considered on a probability space $(\Omega,\mathcal{F},{\sf
P})$, ${\sf E}$ denotes the integration w.r.t. ${\sf P}$.

The main problem is to find a function in genetic and external risk
factors describing the phenotype (that is the individual healthy or
sick) in the best way.

\subsection{Prediction algorithms}\label{pred_alg}

Let $\mathcal X:=\{0,1,2\}^n$ denote the space of all possible
values of explanatory variables.
 Any function $f:\mathcal{X}\to\{-1,1\}$ will be called a {\it theoretical prediction function}.
Define the {\it balanced} or {\it normalized prediction error} for the theoretical prediction function $f$ as
$$Err(f):= \e |Y-f(X)|\psi(Y)$$ where the {\it penalty function} $\psi:\{-1,1\}\to \bbR_+$.
Obviously,
\begin{equation}\label{errorcl0}
Err(f)  =2\psi(-1)\p(f(X)=1,Y=-1)+2\psi(1)\p(f(X)=-1,Y=1).
\end{equation}
Clearly $Err(f)$ depends also on the law of $(X,Y).$ Following
\cite{Velez} and \cite{CrossValid} we put
\begin{equation}\label{psi}
\psi(y)=\frac{1}{4\p(Y=y)},\,\,y\in\{-1,1\},
\end{equation}
the trivial cases $\p(Y=-1)=0$ and $\p(Y=1)=0$ are excluded. Then
\begin{equation}\label{errorcl}
Err(f) =\frac12\p(f(X)=1|Y=-1)+\frac12\p(f(X)=-1|Y=1).
\end{equation} For  a {\it balanced}
sample considered in \cite{MDR}, $\p(Y=-1)=\p(Y=1)=1/2$ and $Err(f)= \e |Y-f(X)|/2$ is equal to the  {\it
classification error} $\p(Y\neq f(X))$.

The reason to consider this weighted scheme is that a
misclassification in a more rare class should be taken into account
with a greater weight. Otherwise, if the probability of disease
$\p(Y=1)$ is small, then the trivial function $f(x)\equiv -1$ may
have the least prediction error. The approach to calculation of
the prediction error based on  penalty functions  is not the only
one possible. Nevertheless Velez et al. \cite{Velez}  showed that
for models with high computational costs it outperforms
substantially other methods such as over- and undersampling.


It is easy to prove  that the {\it  optimal } theoretical prediction
function   minimizing the balanced prediction error is given by
\begin{equation}\label{theo_unb}
f^*(x)=\begin{cases}\;\,1,& p(x)>\p(Y=1),\\
\!\!-1,& \mbox{otherwise.}\\
\end{cases}
\end{equation}
where 
\begin{equation}\label{cprob}
p(x)=\p(Y=1|X=x), \;\;x \in \mathcal X.
\end{equation}
Then each multilocus genotype (with added external risk factors)
$x\in\mathcal{X}$ is classified as high-risk if $f^*(x)=1$ or
low-risk if $f^*(x)=-1.$

Since $p(x)$ and $\p(Y=1)$ are unknown, the immediate application of
\eqref{theo_unb} is not possible. Thus  we try to find an
approximation of unknown function $f^*$ using a {\it prediction
algorithm} that is a function $f_{PA}=f_{PA}(x,\xi(S))$ with values
in $\{-1,1\}$ (recall that $Y\in\{-1,1\}$ a.s.) which depends on
$x\in\mathcal{X}$ and the sample
\begin{equation}\label{samp}
\xi(S)=\{\xi^{j},\;j\in S\}\;\mbox{ where }\; S\subset
\{1,\ldots,N\}.
\end{equation}
 The simplest way is to
employ formula \eqref{theo_unb} with $p(x)$ and $\p(Y=1)$ replaced
by their statistical estimates. For example introduce
\begin{equation}\label{estp}
\widehat{p}(x,\xi(S))= \frac{\sum_{j\in S}I\{Y^j=1,
X^j=x\}}{\sum_{j\in S}I\{X^j=x\}},\;\;x\in \mathcal{X},
\end{equation}
and take
\begin{equation}\label{esty}
\widehat{\p}_S(Y=1)=\frac{1}{\sharp S}\sum_{j\in S}I\{Y^j=1\}
\end{equation}
where $I\{A\}$ stands for the indicator of an event $A$ and $\sharp
D$ denotes the cardinality
  of a finite set $D$.

Along with
\eqref{esty} we will consider
\begin{equation}\label{ecp}
\wh\p_S(Y=1|X\in C)=\frac{\sum_{j\in S}I\{Y^j=1,X^j\in
C\}}{\sum_{j\in S}I\{X^j\in C\}}, \;\;C\subset \mathcal{X}.
\end{equation}
Thus \eqref{estp} is a special case of \eqref{ecp} for $C=\{x\}$
with $x\in\mathcal{X}$.
Note that more
difficult way is to search for the estimators of $f^*$ using several
subsamples of $\xi$.

Assume that we constructed
a prediction algorithm $f_{PA}$. Then taking in mind \eqref{errorcl} set
\begin{equation}\label{pred_alg_error}
Err(f_{PA}(\cdot,\xi(S)))=\frac{1}{2}\sum_{y\in\{-1,1\}}\p\left(f_{PA}(X,\xi(S))\neq
y|Y=y\right).
\end{equation}

As a law of $(X,Y)$ is unknown one can only construct an estimate
$\widehat{Err}(f_{PA}(\cdot,\xi(S)))$ of
$Err(f_{PA}(\cdot,\xi(S)))$. In Section 3 we use the {\it estimated
prediction error} of a prediction algorithm $f_{PA}$ which is based
on $K$-fold cross-validation and has the form
\begin{equation}\label{est_pred_alg_error}
\widehat{Err}_K(f_{PA}(\cdot,\xi),\xi)=\frac{1}{2}\sum_{y\in\{-1,1\}}
\frac1K\sum_{k=1}^K\frac{\sum\limits_{j\in
S_k}I\left\{f_{PA}(X^j,\xi(\overline{S_k}))\neq
y,Y^j=y\right\}}{\sum\limits_{j\in S_k}I\{Y^j=y\}}
\end{equation}
where
\begin{equation}\label{crv}
S_k=\left\{(k-1)\bigg[\frac NK\bigg]+1,\;\ldots,\; k\bigg[\frac
NK\bigg]I\{k<K\}+NI\{k=K\}\right\},
\end{equation}
$\overline{S_k}=\{1,\dots,N\}\setminus S_k$ and $[a]$ is the integer
part of $a\in\mathbb R.$

A very important problem is to make sure that the prediction
algorithm $f_{PA}$ gives statistically reliable results. The quality
of an algorithm is determined by its prediction error
\eqref{pred_alg_error} which is unknown and therefore the  inference
is based on consistent estimates of this error. Clearly the high
quality of an algorithm means that it captures the dependence
between predictors and response variables, so the error is made more
rarely than it would be if these variables were independent.
 Consider a null hypothesis $H_0$ that $X$ and $Y$ are independent.
 If they are in fact   dependent, then
for any significant prediction algorithm $f_{PA}$ an appropriate
test procedure involving $f_{PA}$ should reject $H_0$ at the
approximate  significance level   $\alpha$, e.g., $5\%$.
Intuitively, this shows that results of the algorithm could not be
obtained by chance. For such a procedure, we take a {\it permutation
test } (see \cite{Golland05}).  Its idea is as follows.

Permutation test for a given statistic $\wh L(\xi)$ (we consider $\wh L(\xi)=\widehat{Err}(f_{PA}(\cdot,\xi))$
is done by the following
steps.

\begin{enumerate}
\item Generate $B$ independent random vectors
$(\pi_1^b,\ldots,\pi_N^b),$ $1\le b\le B,$ with the uniform
distribution over all permutations $\Pi_N$ of $1,\ldots,N$.
\item
Compute $\widehat{Err}_{K,b}=\widehat{Err}_K (f_{PA},\bar\xi_b)$,
$1\le b\le B,$ with
$$\bar\xi_b=\big(\big(X^{1},Y^{\pi_1^b}\big),\ldots,\big(X^{N},Y^{\pi_N^b}\big)\big).$$
%
\item Find the {\it Monte Carlo ${\sf p}$-value} (see, e.g., \cite[p. 63]{Lehmann}):
\begin{equation}\label{montecarlopvalue}
\widehat {\sf p}=\widehat{ F}\big( \widehat{Err}_K(f_{PA},\xi)\big)
\end{equation} where $\wh F=\wh F(z)$ is the empirical cumulative distribution function (c.d.f.)
defined by the relation
$$\wh F(z)=\frac1B\sum_{b=1}^BI\big\{\widehat{Err}_{K,b}\le z\big\},\;\;
z\in\mathbb R.$$
\item If $\wh {\sf p}<\alpha$, reject $H_0,$ otherwise not.
%
\end{enumerate}

According to \cite{Golland05}, one ideally has to use all
permutations belonging to $\Pi_N$ but this is impractical  in view
of computational costs. Thus the Monte Carlo approximations for the
{\it true} ${\sf p}$-{\it value} ${\sf
p}=F\big(\widehat{Err}_K(f_{PA},\xi)\big)$ are employed, here $F$ is
the c.d.f. of $\wh {Err}_{K,b}$. The  upper
bound for $|{\sf p}-\widehat{\sf p}|$ is $1/2\sqrt{B}$ (see
\cite{Golland05}). This could be used to determine the number $B$ of
simulations for a desired accuracy.

Note also that if the estimate  of the error function for the
algorithm $f_{PA}(\cdot,\xi)$ is {\it asymptotically optimal}, i.e.
converges in probability to the error of the optimal prediction
function $f^*$ as $N\to\infty$ ($\xi$ depends on $N$), then the rule
of thumb is  to suspect overfitting if $\widehat{Err}_K(f_{PA},\xi)$
is close to $1/2,$ which is a probability limit of this error under
$H_0$ as $N\to \infty$.

We use complementary approaches to analyze  dataset related to
complex diseases. Each approach (MDR, LR and machine learning) is
characterized by its own way of constructing prediction algorithms.
For each method one or several prediction algorithms admitting the
least estimated prediction error are found (a typical situation is
that there are several ones with almost the same estimated
prediction error). These prediction algorithms provide a way to
determine the domains where the disease risk is high or low
(depending on the value of the corresponding prediction function).
It is also possible to select combinations of SNPs and external risk
factors whose presence influences the liability to disease to a
great extent. Some methods allow to present such combinations
immediately. Others, which employ more complicated forms of
dependence between explanatory and response variables,  need further
analysis based on modifications of permutation tests.

Now we pass to the description of various statistical methods and
their applications to the cardiovascular risk detection.

\subsection{Multifactor dimensionality reduction}

Ritchie et al. \cite{MDR} introduced {\it multifactor dimensionality
reduction} (MDR) as a new method of analyzing gene-gene and
gene-environment interactions. Rather soon the method  became very
popular. Since the first publication   more than 200 papers applying
MDR in genetic studies were written   (see, e.g., references in
\cite{Velez}).

MDR is a flexible non-parametric method  not depending on a
particular
inheritance model. 
We give a rigorous description of the method following ideas of
\cite{MDR} and \cite{Velez}. As mentioned earlier, the probability
$p(x)$ introduced in \eqref{cprob} is unknown. To find its estimate
one can apply maximum likelihood approach assuming that the random
variable $I\{Y=1\}$ conditionally on $X=x$ has a Bernoulli
distribution with unknown parameter $p(x)$. Then we come to
\eqref{estp}.

A direct calculation of estimate in
\eqref{estp} with exhaustive search over all possible values of $x$
is  highly inefficient, since the number of different values of $x$
grows   exponentially with number of risk factors. Moreover, such a
search leads to overfitting. Instead, it is usually   supposed that
 $p(x)$ non-trivially depends    not on all, but certain
 variables $x_i$. That is, there exist $l\in\mathbb{N},\,l<n,$ and $(k_1^*,\dots,k_l^*),\,$ where $1\le k^*_1<\ldots<k^*_l \le n$, such that
   for each $x=(x_1,\ldots,x_n) \in
\mathcal{X},$ the following relation holds:
\begin{equation}\label{part}
p(x)=\p(Y=1|X_{k^*_1}=x_{k^*_1},...,X_{k^*_l}=x_{k^*_l}).
\end{equation}
 In other words only few factors influence the disease and the
others can be neglected. A minimal combination of factors
$(X_{k^*_1},\dots,X_{k^*_l})$ in formula \eqref{part} is called {\it
the most significant}. Clearly it is the most significant
combination which has the least prediction error. Indeed, if we
consider any other combination of pairwise different indices
$k_1,\dots,k_r$ and set
\begin{equation*}f_{k_1,\ldots,
k_r}(x)=\begin{cases}\;\;1, &
\p(Y=1|X_{k_1}=x_{k_1},...,X_{k_r}=x_{k_r})>\p(Y=1),\\
\!-1, & \mbox{otherwise,}\end{cases}
\end{equation*}
then we obviously have
\begin{equation}\label{other_comb}
Err\left(f_{k^*_1,\ldots, k^*_l}\right)\le Err(f_{k_1,\ldots, k_r})
\end{equation}
where $Err(f)$ is calculated according to \eqref{errorcl}.

To choose the most significant combination, exhaustive search over
all possible combinations of factors is applied. 
For each $\{k_1,\dots,k_r\}\subset\{1,\ldots,n\}$ and any
$x\in\mathcal{X}$ consider
$$C_{k_1,\dots,k_r}(x)=\{u=(u_1,\ldots,u_n)\in\mathcal{X}:u_{k_i}=x_{k_i},\,i=1,\dots,r \}$$
and for $S$ appearing in \eqref{samp} define a prediction algorithm
(cf. \eqref{theo_unb}) by
\begin{equation}\label{pred_MDR}
\widehat{f}_{k_1,\ldots, k_r}(x,\xi(S)):=\begin{cases} \;\;1, &
\widehat{\p}_S(Y=1|X\in C_{k_1,\dots,k_r}(x))>\widehat{\p}_S(Y=1),\\
\!-1,& \mbox{otherwise,}
\end{cases}\end{equation}
here we use formulas \eqref{esty} and \eqref{ecp}.
 It is easy to show that $\widehat{Err}_K(\widehat{f}_{k_1,\ldots, k_r},\xi)\to Err(f_{k_1,\ldots, k_r})$ in probability as
$N\to\infty$ ($\xi$ depends on $N$). Consequently, \eqref{other_comb} implies that, for any $\varepsilon>0$ and all $N$ large
enough, with
probability close to $1$ one has
\begin{equation*}
\widehat{Err}_K(\widehat{f}_{k^*_1, \ldots,
k^*_l},\xi)<\widehat{Err}_K(\widehat{f}_{k_1,\ldots, k_r},\xi)+\varepsilon.
\end{equation*}
Hence,  it is natural to pick one or a few combinations of factors
with the smallest empirical prediction errors 
as an approximation for the most significant combination. 


The last step in MDR  is to determine statistical significance of
the results. Here we test a null hypothesis of independence between
$X$ and $Y$ i.e. between risk factors $X$ and a disease $Y$. This
can be done via the permutation test described in Section 2.1.

\vskip0,2cm \textbf{\textbf{MDR method with ``independent rule''}}.
We propose {\it multifactor dimensionality reduction} with {\it
``independent rule''} (MDRIR) method to improve the estimate of
probability $p(x)$. This approach is motivated by Park
\cite{Park09}, who deals with classification of large array of
binary data. The principal difficulty with employment of  formula
\eqref{estp}
 is that the number of observations in numerator and
denominator of the formula might be small even for large $N$ (see,
e.g., \cite{Lee}). This can lead to inaccurate estimates and finally
to a wrong prediction algorithm. Moreover, for some samples the
denominator of
\eqref{estp} can equal zero.

The Bayes formula implies that
\begin{equation}\label{Indrule1}
p(x)=\frac{\p(X=x|Y=1)\p(Y=1)}{\p(X=x|Y=1)\p(Y=1)+\p(X=x|Y=-1)\p(Y=-1)},
\end{equation}
where the trivial cases $\p(Y=-1)=0$ and $\p(Y=1)=0$ are excluded.
Substituting \eqref{Indrule1} into \eqref{theo_unb} we obtain the
following expression for prediction function:
\begin{align}\label{solution_ind_MDR}
f^*(x)= \begin{cases} 1,&\p(X=x|Y=1)>\p(X=x|Y=-1),\\
-1,& \mbox{otherwise}. \end{cases}\end{align}

As in standard MDR method described above, we will assume that
formula \eqref{part} holds. It was proved in \cite{Park09} that for
a broad class of models (e.g., {\it Bahadur model} \cite{Bahadur61},
{\it logit model} \cite{Logit}) the conditional probability
$\p\left(X_{k_1}=x_1,\dots, X_{k_r}=x_r\left|Y=y\right.\right)$,
where $y=\pm1$, can be estimated in the following way:
\begin{equation}\label{Indrule2}
\widehat{\p}_S\left(X_{k_1}=x_1,\dots,
X_{k_r}=x_r\left|Y=y\right.\right)=\prod\limits_{i=1}^r
\widehat{\p}_S\left(X_{k_i}=x_i\left|Y=y\right.\right),
\end{equation}
here (cf. \eqref{ecp})
\begin{equation}\label{Indrule2a}
\widehat{\p}_S\left(X_{k_i}=x\left|Y=y\right.\right)=\frac{\sum\limits_{j\in
S} I\{X_{k_i}^j=x, Y^j =y\}}{\sum\limits_{j\in S} I\{Y^j=y\}}.
\end{equation}

Combining \eqref{part}, \eqref{solution_ind_MDR}, \eqref{Indrule2} and
\eqref{Indrule2a} we find the desired estimate of $f^*(x)$.

A number of observations in numerator and denominator of
\eqref{Indrule2a}  increases considerably comparing with
\eqref{pred_MDR}. It allows to estimate the conditional probability
more precisely whenever the estimate introduced in \eqref{Indrule2}
is reasonable.

Thus, as opposed to standard MDR method, MDRIR uses alternative
estimates of conditional probabilities. All other steps (prediction
algorithm construction,  prediction error calculation) remain the
same. Let us mention that as far as we know
this modification of MDR  has not been applied before. It is based
on a combination of the original MDR method \cite{MDR} and the ideas
of \cite{Park09}.

\subsection{Logic regression}

The {\it logic regression} (LR)  was proposed in \cite{RuczKB}.
Further generalizations are given in \cite{FriIckst},
\cite{KooperBis}, \cite{Schwender2}, \cite{Schwender}
  and other works. LR is based on the
classical binary logistic regression \cite{logreg} and exhaustive
search for relevant predictor combinations. The main difficulty is
to organize the search  efficiently. The LR method was applied to
identification of the most significant SNP combinations in
 \cite{Albrechtsen}, \cite{Nunke} and \cite{Schwender2}. Note that for genetic analysis it is convenient to use explanatory variables taking
 3 values. Thus we employ {\it ternary
 variables}, whereas the authors of the above-mentioned papers employ   binary ones.

Let $p(x)$ be  the conditional probability  of a disease defined in
\eqref{cprob}. We suppose that trivial situations when
$p(x)\in\{0,1\}$ do not occur and omit them from the consideration.
To estimate $p(x)$  we pass now to the {\it logistic transform}
\begin{equation}\label{logistpr}q(x)=\lambda\left(p(x)\right)\end{equation}
where $\lambda(z)=\ln(z/(1-z)),\,z\in(0,1),$ is the {\it inverse
logistic function}. The {\it logistic function} itself equals to
$\Lambda(t)=(1+e^{-t})^{-1},\,t\in\bbR.$ Note that we are going to
estimate the unknown disease probability with the help on linear
statistics with appropriately selected coefficients. Therefore it is
natural to avoid restrictions on possible values of the function
estimated. Thus the logistic transform is convenient, because
$p(x)\in(0,1)$ for $x\in \mathcal{X}$ while $q(x)$ can take all real values.

Consider a class $\mathcal{G}$  of all real-valued functions in
ternary variables $x_1,\dots,x_n$. We call a {\it model} of the
dependence between the  disease and explanatory variables any
subclass $\mathcal{M}\subset
\mathcal{G}$. Set
$$\wh\psi(y,\xi(S))=\frac{1}{4\widehat{\p}_S(Y=y)},\,\,\,y\in\{-1,1\},
$$ here $\widehat{\p}_S(Y=y)$ was introduced in \eqref{esty}. Define the {\it normalized  smoothed score function}
\begin{equation}\label{sglazh}L(h,\xi(S))= \frac{1}{\sharp S}\sum_{j\in S}\phi(-Y^jh(X^j))\wh\psi(Y^j,\xi(S))
  \end{equation} where $S$ is introduced in \eqref{samp}, $\phi(t)=\log_2(1+e^{t})$ for $t\in\mathbb{R},$
 and $h\in
\mathcal{M}$. In contrast to previous works our version of LR scheme
involves normalization (cf. \eqref{errorcl0}), i.e. taking the
observations with weights dependent on the proportion of
cases and controls in subsample $\xi(S)$.

An easy computation yields that  $\arg\min_{h\in \mathcal{M}}L(h,\xi(S))$ equals
to
$$\arg\max_{h\in \mathcal{M}}\frac{1}{\sharp S}\sum_{j\in S}\left(\ln\Lambda(h(X^j))
\frac{I\{Y^j=1\}}{2\widehat{P}_S(Y=1)}+
\ln(1-\Lambda(h(X^j)))\frac{I\{Y^j=-1\}}{2\widehat{P}_S(Y=-1)}\right).
$$
That is, minimizing the score function is equivalent to the search
of normalized maximal likelihood estimate of $q$. Note that estimating the
disease probability in this setup is closely connected with the
problem of data classification, i.e. predicting the disease by the
value of  $x\in \mathcal{X}.$ Recall that in standard
 classification problem instead of the score function
\eqref{sglazh} one uses the following normalized estimate of the error
probability
$$\widetilde{L}(h,\xi(S))=\frac{1}{\sharp S}\sum_{j\in S} I\{Y^jh(X^j)<0\}
\wh\psi(Y^j,\xi(S)). $$ 
In fact the optimal choice of $h$ for these problems coincide if the
underlying model $\mathcal{M}$ is correctly specified (i.e. $q\in \mathcal{M}$), see
\cite{Biau08}. However the usage of score function $L$ has an
important advantage over $\widetilde{L}$ because one has to evaluate
the minimum of a smooth function.

\begin{figure}
\begin{center}
\includegraphics[width=4.7cm]{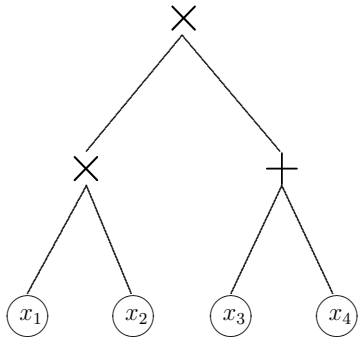}

\caption{A tree $T$ representing a function $T(x)=(x_1\times
x_2)\times(x_3+x_4).$} \label{onetree}
\end{center}
\end{figure}

A wide and easy to handle class of models is obtained by taking
functions linear in variables $x_1,\dots,x_n$ or in their products.
In turn these functions admit a convenient representation by
elementary polynomials. Recall that an {\it elementary polynomial } (EP) is a
function  $T$ in ternary variables  $x_1,\dots,x_n$ belonging to $
\{0,1,2\}$ which can be represented as a finite sum of products
$x_1^{u_1}\dots x_n^{u_n}$ where $u_1,\dots,u_n\in\mathbb{Z}_+.$
The addition and multiplication of ternary variables is considered
by modulo 3. Any EP can be represented as a  {\it binary tree}\footnote{For the basic concepts of the graph theory see, e.g., \cite{Bondy}.} in
which {\it knots} (vertices which are not {\it leaves}) contain either
addition or multiplication sign, and each leaf corresponds to a
variable. Figure \ref{onetree} provides an example of a binary
tree. Different trees may correspond to the same EP, thus this
relation is not one-to-one. However, it does not influence our
problem, so we regain the notation $T$ for a tree. A finite set of
trees $F=(T_1,\dots,T_s)$ is called a {\it forest}. For a tree $T,$
its {\it complexity } $C(T)$ is the number of leaves. The complexity
$C(F)$  of a forest $F$ is the maximal complexity of trees
constituting $F.$

It is clear that if $g\in\mathcal{G}$ then there exists $s\ge1$ such
that $g$ has the following form:

\begin{equation}\label{logic_regr_func}
g(x_1,\dots,x_n)=\beta_0+\sum_{i=1}^s \beta_iT_i(x_1,\dots,x_n),
\end{equation}
here $\beta_0,\beta_1,\dots,\beta_s\in\mathbb{R}$ and $T_1,\dots,T_s$ are EP.

Let us say that function $g$ belongs to a class $\mathcal{G}_r(s)$, where
$s,r\in\mathbb{N}$, if there exist a decomposition
\eqref{logic_regr_func} of $g$ such that all trees $T_i$ ($i=1,\ldots,s$) have
complexity less or equal $r$. We identify  a function $g\in
\mathcal{G}_r(s)$ with pair $(F,\beta)$ where $F$ is the
corresponding forest and $\beta=(\beta_0,\dots,\beta_s)$ is the
vector of coefficients in  \eqref{logic_regr_func}.

Minimization of $L(h,\xi(S))$ defined by \eqref{sglazh} over all
functions $h\in \mathcal{M}\subset \mathcal{G}_r(s)$ is done in two
alternating steps. First we find the optimal value of $\beta$ while
$F$ is fixed (which is the minimization of a smooth function in
several variables) and then we search for the best $F$. Here one
uses stochastic algorithms, since the number of such forests
increase rapidly when the complexity $r$ grows. For $s\in\mathbb{N}$,
 a forest $F=(T_1,\dots,T_s)$ and a subsample $\xi(S)$ (see
\eqref{samp}) consider a prediction algorithm $f^F_{\rm{LR}}$
setting
\begin{equation*}\label{lr_predalg} f^F_{\rm{LR}}(x,\xi(S))=\begin{cases} \;\,1,&
\widehat{h}(x)>0, \\\!\! -1,& \mbox{ otherwise,}\end{cases}
\end{equation*} where $\widehat{h}=(F,\widehat{\beta})$ and
\begin{equation}\label{lr_predalg2}\widehat{\beta}=
 \arg\min_{\beta}L\left(\beta_0+\sum_{j=1}^s \beta_j T_j(\cdot),\xi(S)\right).\end{equation} 

Define also the {\it normalized prediction error of a forest}
$F=(T_1,\dots,T_s)$ as
$$\widetilde{\varphi}(F)=\widehat{Err}_K(f^F_{\rm{LR}}(\cdot,\xi),\xi).  $$

A subgraph $B$ of a tree $T$ is called a {\it branch} if it is
itself a binary tree (i.e. it can be obtained by selecting one
vertex of   $T$ together with its offspring). 
Sum and product signs standing in a knot of a tree are called {\it
operations}, thus $\ast$ stands for sum or product.
 Following   \cite{RuczKB}, call  the tree $\widetilde{T} $   {\it a
neighbor } of $T$ if it is obtained  from $T$ via one and only one
of the following transformations.

\begin{enumerate} \item Changing one variable to another in a leaf of the tree $T$ ({\it{variable change}}).

\item Replacing an operation in a knot of a tree $T$ with another   one, i.e. sum to product or vice versa ({\it{operator change}}).

\item Changing  a branch of two leaves to one of these leaves ({\it{deleting a leaf}}).

\item Changing a leaf to a branch of two leaves, one of which contains the same variable as in initial leaf ({\it{splitting a leaf}}).

\item Replacing  a branch $B_1\ast B_2$ with the branch  $B_1$ ({\it{branch
pruning}}).

\item Changing a branch  $B$ to a branch  $x_j\ast B$ ({\it{branch growing}}), here  $x_j$ is a variable.

\end{enumerate}

\begin{figure}
\begin{center}

\includegraphics[width=10cm]{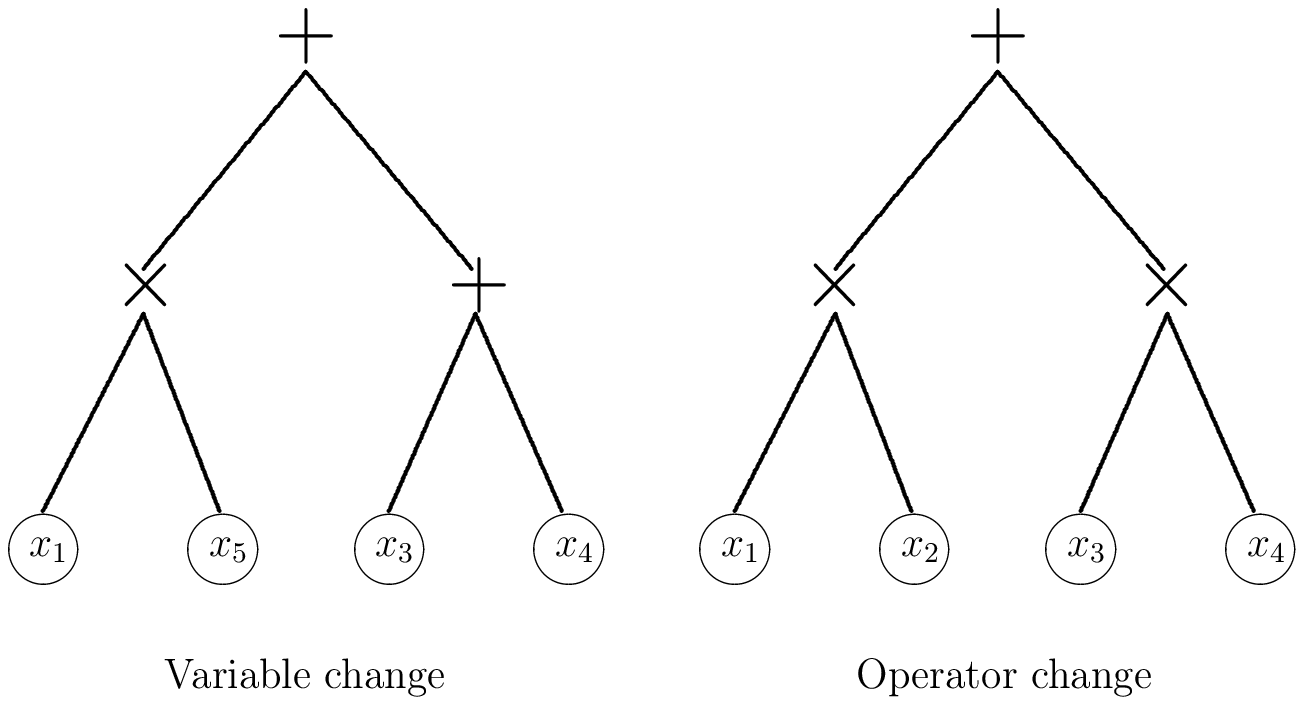}
\vspace{5mm}
\includegraphics[width=10cm]{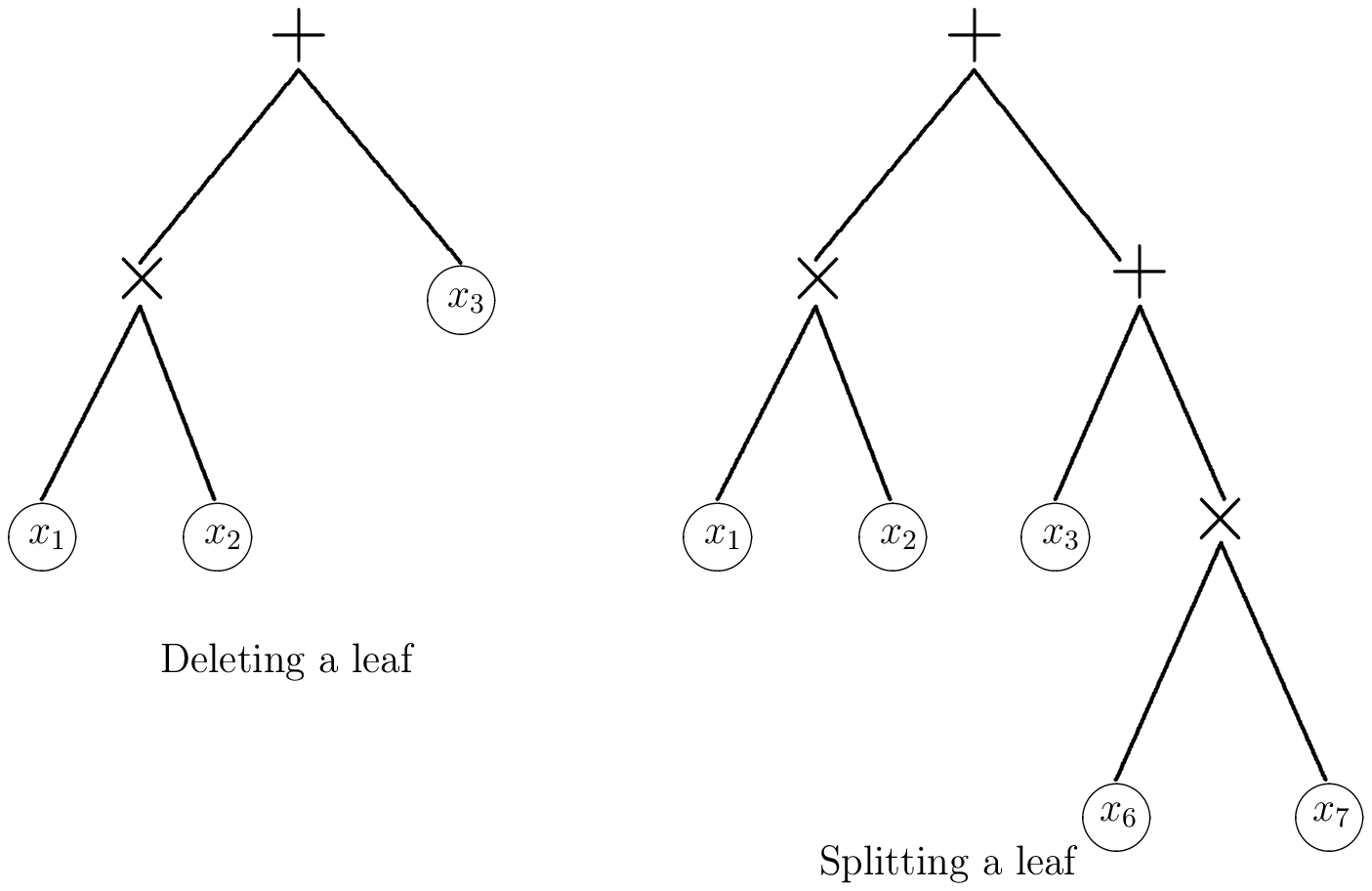}
\vspace{5mm}
\includegraphics[width=10cm]{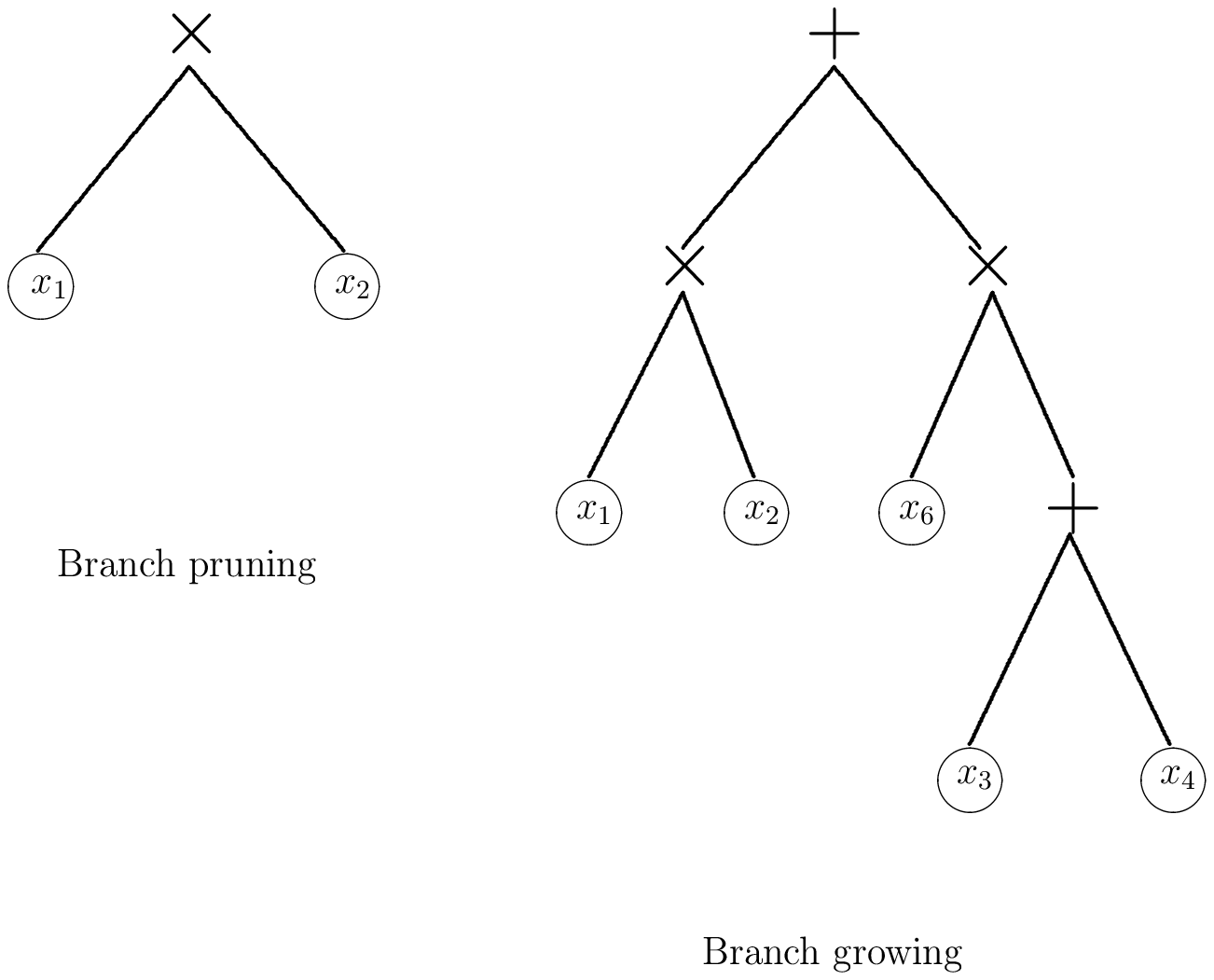}
\caption{Neighbors of the tree $T.$}\label{multree}
\end{center}
\end{figure}

Figure \ref{multree} depicts results of these operations applied to
the tree $T$ of Figure \ref{onetree}. We say that forests $F$ and
$\widetilde{F}$ are  {\it neighbors} if they can be written as
$F=\{T_1,T_2,\dots,T_s\}$ and
$\widetilde{F}=\{\widetilde{T_1},T_2,\dots,T_{s}\}$ where  $T_1$
and $\widetilde{T_1}$ are neighbors. The neighborhood relation
defines a finite connected graph on all forests of equal size $s$
with complexity not exceeding  $r$. To each vertex $F$ of this graph
we assign a number $\widetilde{\varphi}(F)$. To
 find the global minimum of a function defined  on a finite
 graph we employ the   {\it simulated annealing method}
 (see, e.g., \cite{Hajek}, \cite{Nikolaev} and \cite{Sylvain09}).
 This method  constructs some specified Markov
  process which takes values in the graph vertices and converges with high probability
  to the global minimum
of the function.
 To avoid stalling at a local minimal point the process  is allowed
 to pass with some small probability to a point   $F$ having
 greater value of   $\widetilde{\varphi}(F)$ than current  one.
 We propose a new modification of this method
 in which the output is the forest corresponding to the minimal
 value of a function $\widetilde{\varphi}(F)$ over all (randomly) visited points.

\subsection{Machine learning methods}

Let us describe (see, e.g., \cite{MachLearn09}) two among the most
popular machine learning methods -- {\it random forests} (RF) and
{\it stochastic gradient boosting} (SGB). They belong to {\it
ensemble methods} which combine multiple predictions from a certain
base algorithm to obtain better predictive power (i.e. less
estimated prediction error). We use {\it classification and regression trees}
(CART) for a base learning algorithm because it showed good
performance in a number of studies (see \cite{Hastie10}).

{\it Classification tree} $T$ is a binary tree having the following
structure. Any leaf of $T$ contains either $1$ or $-1$ and for any
vertex $P$ in $T$ (including leaves) there exists a subset $A_P$ of
the explanatory variable space $\mathcal X$ such that the following
properties hold:

\begin{enumerate} \item
$A_P=\mathcal X$ if $P$ is the root of $T$;

\item if vertices $P'$ and $P''$ are children for $P$, then $A_{P'}\cup
A_{P''}=A_P$ and $A_{P'}\cap A_{P''}=\varnothing.$
\end{enumerate}

\begin{figure}
\begin{center}
\includegraphics[width=6cm]{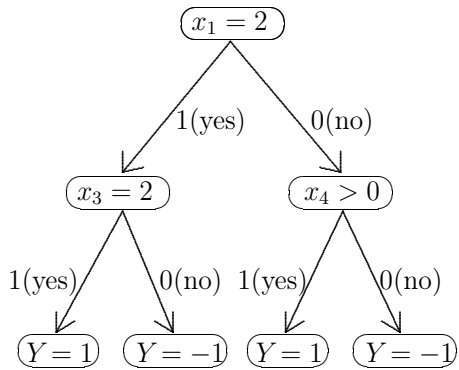}
\caption{CART representing the prediction $((x_1=2)$ and $(x_3=2))$
or $(x_4>0).$} \label{Img:carttree}
\end{center}
\end{figure}

In particular, subsets corresponding to the leaves form the
partition of $\mathcal X$. To obtain a prediction of $Y$ given a
certain value $x\in \mathcal{X}$ of the random vector $X$, one
should go along the path which starts from the root and ends in some
leaf turning at each parent vertex $P$ to that child $P'$ for which
$A_{P'}$ contains $x$. At the end of the $x$-specific path, one gets
either $1$ or $-1$ which serves as a prediction of $Y$. Figure \ref{Img:carttree}
provides an example
of a classification tree.
Namely, the partition of $\mathcal{X}$ is formed by values of
boolean functions standing in parent vertices. For each $x$ starting
from the root of the tree we calculate the value of a boolean
function and move along the edge marked with the value obtained ($1$
or $0$). The left child of the root corresponds to the subset
$\{x\in\mathcal{X}:x_1=2\},$ while the right one to its complement
in $\mathcal{X}.$ Next, the leftmost leaf stands for a subset
$\{x\in\mathcal{X}:x_1=2,x_3=2\},$ and if $X$ falls in this subset,
we predict that $Y=1;$ the rightmost leaf stands for a subset
$\{x\in\mathcal{X}:x_1<2,x_4=0\},$ and if $X$ takes values in this
subset,   we predict  $Y=-1.$

Classification tree could be constructed via CART {\it algorithm} (if it
is the case, we will call it CART). The algorithm proceeds
iteratively. That is, on the $l$-th step of the algorithm
($l=1,2,\ldots$), each element $A$ of the current partition
$\mathcal A_l$ ($\mathcal A_1=\mathcal X$) of the set
$\mathcal X$ is divided into two disjoint parts
$$A^+(i,t)=\{(x_1,\dots,x_n)\in A:x_{i}\leq t\}\mbox{ and } A^-(i,t)=\{(x_1,\dots,x_n)\in A:x_{i}>t\}$$
minimizing the sum $ \wh G(A^+(i,t))+\wh G (A^-(i,t))$ over
$i=1,\ldots,n$ and $t\in\{0,1\}$. Here the {\it empirical Gini
index}
$$\wh G(C)=2\wh\p_S(Y=1|X\in C)\bigl(1-\wh\p_S(Y=1|X\in C)\bigr)$$
with $C\subset\mathcal X$ and
$\wh\p_S(Y=1|X\in C)$ (see \eqref{ecp}) measures the heterogeneity
of the subsample $\{j\in S:X^j\in C\}$ w.r.t.
response variable
$Y$. Any uninformative partition with
$$\min_{(i,t)}\bigl(\wh G(A^+(i,t))+\wh G(A^-(i,t))\bigr)> \wh G(A),$$ is not made.

The algorithm stops whenever a number of leaves $D$ reaches some
critical level  which is chosen via some data-dependent criteria
(see \cite{Hastie10}, page 308). For a subsample $\xi(S)$ of $\xi,$
each CART defines a prediction algorithm
\begin{equation}\label{cart_predalg} f(x,\xi(S))=\sum_{d=1}^D a_d(\xi(S))
I\{x\in A_d(\xi(S))\}\end{equation} where
$\{A_1(\xi(S)),\dots,A_D(\xi(S))\}$ is the partition of $\mathcal X$
corresponding to the leaves,
\begin{equation*}
a_d(\xi(S))=\left\{%
\begin{array}{ll}
1, & \sharp\{j\!\!\in \!S\!:Y^j=1,X^j\!\in\!
A_d(\xi(S))\}>\sharp\{j\!\!\in\! S\!:Y^j=-1,X^j\!\in\! A_d(\xi(S))\}; \\
\!\!\!\!-1, & \mbox{otherwise}. \\
\end{array}%
\right.
\end{equation*}

RF is a non-parametric method of estimating conditional probability
$p(x)$. It was successfully applied to genetics data in a number of
papers (see references in \cite{MachLearn09}).
It could be  briefly described as follows (see chapter 15 in
\cite{Hastie10} for details). Generate  $B$ bootstrap samples from
the initial sample where one could choose $B=\max\{[N\log N],1000\}$
according to \cite{Olive}. For $b$-th bootstrap sample $(1\le b\le
B)$ construct a CART prediction algorithm   $f_b:\mathcal
X\times(\mathcal X\times\{-1,1\})^N\to\{-1,1\}$ defined according to
\eqref{cart_predalg} and take
$$\wh
p_{\rm{RF}}(x,\xi(S))=\big(B^{-1}\sum_{b=1}^Bf_b(x,\xi(S))+1\big)/2$$
as an estimate of $p(x)$.

It is shown in \cite{Biau08} that generally RF method gives
consistent estimates of  $p(x)$  only if the  number  of partitions
used in CART grows slower than the sample size.
A final
prediction algorithm $f_{\rm{RF}}(x,\xi(S))$ is constructed from the
estimate $\wh p_{\rm{RF}}(x,\xi(S))$ similarly to \eqref{theo_unb},
i.e.
\begin{equation*}f_{\rm{RF}}(x,\xi(S))=\begin{cases}\;\,1, &\wh p_{\rm{RF}}(x,\xi(S))>\wh \p_S(Y=1),\\\!\! -1,&
\mbox{ otherwise}.\end{cases}\end{equation*} The distinctive
features of this method are low computational costs and the ability
to extract relevant predictors when the number of irrelevant ones is
large (see \cite{Biau10}).

SGB is another non-parametric method of estimating conditional
probability  $p(x)$. This method is used in a number of procedures
for studying genetics
data (see, e.g., \cite{SNPHunter09}). SGB  method can be described as follows (\cite{Friedman01}).\\

\begin{enumerate}
\item Pass on the input of the algorithm\footnote{This algorithm is not to be confused with prediction algorithms.}
initial parameters $D,M\in\mathbb N$ and $\rho,\eta\in(0,1)$.\\

\item Put $m=0$, $\xi_0(S)=\xi(S)$ and  $$f_0(x,\xi_0(S))\equiv
\frac12\ln\frac{\wh\p_S(Y=1)}{\wh\p_S(Y=-1)}.$$

\item Increase $m$ by 1 and define
$$\bar Y^j_m:=\frac{2Y^j}{1+\exp\{2Y^jf_{m-1}(X^j,\{\xi_l(S)\}_{l=0}^{m-1})\}}.$$
\vskip0,2cm Choose a random subset in $\xi_m(S)=\{(X^j,\bar
Y_m^j)\}_{j\in S}$ with $[\eta \sharp S]$ elements.  Construct CART
prediction algorithm (with $D$ leaves)  $\sum_{d=1}^Da_d^m(\xi_m(S))
I\{x\in A_d^m(\xi_m(S))\}$ on the chosen subset. Compute weight
coefficients
$$w_d^m(\xi_m(S))=\frac{\sum_{j\in J   }\bar Y^j_m}{\sum_{j\in J }|\bar Y^j_m|(2-|\bar Y^j_m|)},\quad d=1,\dots,D,$$
where the random set $J=\{j:X^j\in A_d^m(\xi_m(S))\},$ and put
$$f_m=f_{m-1}+\rho\sum_{d=1}^D w_d^m(\xi_m(S))I_{A_d^m(\xi_m(S))},$$
here $\rho$ is the {\it memory relaxation} parameter.

\item If $m<M$, go to Step 3, otherwise determine a final estimate
$$\wh p_{\rm{SGB}}(x,\xi(S))=\frac1{1+\exp\{-2f_M(x,\xi(S),\{\xi_m(S)\}_{m=1}^M)\}}.$$

\end{enumerate}

This algorithm is to be run for several times with different
parameters $D$, $M$, $\rho$ and $\eta$. Then their optimal values
could be chosen via cross-validation (see section 16.3.1 in
\cite{Hastie10}). Small values of $\eta$ ($=0.1,\,0.05,\,0.0225$
etc.)  help to get accurate estimates for relatively noisy data.

Standard RF and SGB work poorly for unbalanced samples. One needs
either to balance given datasets (as in \cite{Chawla10}) before
these methods are applied or use special modifications of RF
(\cite{Breiman04}) and SGB (\cite{Liu08}). To avoid overfitting,
permutation test needs to be done.

A common problem of all machine learning methods is a complicated
functional form of the final probability estimate  $\wh p(x,\xi)$ (w.r.t. $x$). In
genetic studies, one wants to pick up all relevant combinations of
SNP and risk factors, based on a biological pathway causing the
disease. Therefore, the final estimate $\wh p(x,\xi)$  is to be
analyzed.
\vskip0,3cm
 We describe  one of the possible methods of such analysis  within RF framework and
called {\it{conditional variable importance measure}} (CVIM). One
could determine CVIM for each predictor $X_i$ in $X$ and range all
$X_i$ in terms of this measure. Following \cite{Strobl08}, CVIM of
predictor $X_i$ given certain subvector $Z_i$ of $X$ is calculated
as follows (supposing $Z_i$ takes values $z_{i1},\ldots,z_{im(i)}$).
\begin{enumerate}
\item Construct a vector $(l_1,\dots,l_N)$, randomly permuting $1,\dots,N$ in each subset
$$A_{ik}=\{j:Z_i^j=z_{ik}\},\quad k=1,\ldots, m(i).$$
\item Generate $B$ bootstrap samples
$\xi_b=\big((X^{jb},Y^{jb}),j=1,\ldots,N\big),$ $b=1,\ldots,B.$ For
each of these samples, construct a classifier $f_b(x,\xi_b)$  and
calculate
\vskip0,2cm
$${\text{CVIM}}_b=\frac1{|C_b|}\sum_{j\in C_b}I\{Y^j=f_b(X^j,\xi_b)\}-
 \frac1{|C_b|}\sum_{j\in C_b}I\{Y^j=f_b(X^{l_j},\xi_b)\}$$
where  $C_b=\{j\in\{1,\ldots,N\}:(X^j,Y^j)\notin \xi_b\}$.
\vskip0,2cm
\item Compute the final CVIM using the formula
\begin{equation}\label{cvim}
CVIM=B^{-1}\sum_{b=1}^B {\text{CVIM}}_b.
\end{equation}
\end{enumerate}
\vskip0,2cm Any permutation $(l_1,\dots,l_N)$ in the CVIM algorithm
destroys dependence between $X_i$ and $(Y,Z_{-i})$ where $Z_{-i}$
consists of all components of $X$ which are not in $Z_i$. At the
same time it preserves initial empirical distribution of $(X_i,Z_i)$
calculated for the sample $\xi$. After that the average loss of
correctly classified $Y$ is calculated. If it is relatively large
w.r.t. CVIM of other predictors, then $X_i$ plays important role in
classification and vice versa.

For $Z_i$, one could take all components $X_k$ ($k\neq i$) such that
the hypothesis of the independence between $X_k$ and $X_i$ is not
rejected at some significance level (e.g., $5\%$). Note also that
CVIM-like algorithm could be used to range pairs of SNP and risk
factors w.r.t. the level of association to the disease. This will be
done elsewhere.

\vskip1,5cm

\section{Applications: risks of IHD and MI}

We employ here the various statistical methods described above to
analyze the influence of genetic and external factors on risks of
ischemic heart disease (IHD) and myocardial infarction (MI) using
the data for  454 individuals (333 cases, 121 controls) and 333
individuals (165 cases, 168 controls) respectively. These data
contain   values of seven SNPs (PAI-1, GpIa, GpIIIa, FXIII, FVII,
IL-6, Cx37), as well as four external risk factors, namely, obesity
(Ob), arterial hypertension (AH), smoking (Sm) and
hypercholesterolemia (HC). The age of all individuals in case and
control groups ranges  from 35 to 55 years, which reduces its
influence   on the risk analysis.
 For each of considered methods, $K$-fold cross-validation is used with $K=6$. As shown in
  \cite{MotRit2006} and \cite{CrossValid}, the standard choice of partition
   number of cross-validation from 6 to 10
 does not change the prediction error significantly. We take $K=6$
 as the sample sizes do not exceed 500.
 The supercomputer SKIF MSU  ``Chebyshev'' was involved to perform computations.
 All applied methods have prediction error    less than 0.25,
 so predictions  constructed have significant predictive power. Indeed,
 in \cite{Coffey} and
\cite{ChiaTiTsai} the interplay between genotype characteristics and
MI development was also studied, with estimated prediction errors
0.30--0.40. Further on we write prediction error instead of estimated prediction error.

\vspace{5mm}

\subsection{MDR and MDRIR methods}

{\bf Ischemic heart disease} \vspace{2mm}

\begin{figure}[!h]

   \begin{center}

\includegraphics[width=8cm]{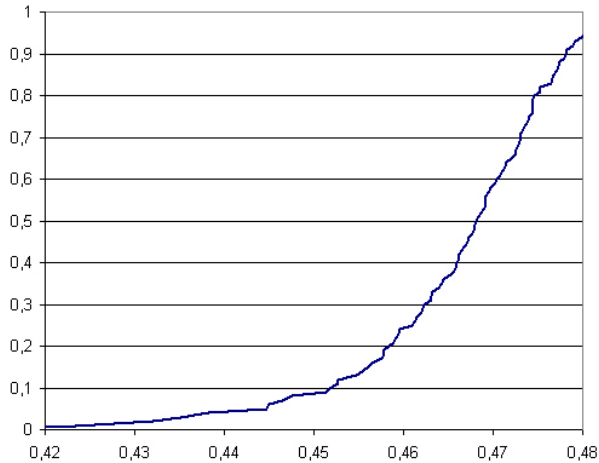}
   \end{center}
\caption{The empirical distribution function of the estimated error for the case when disease risk and
predictors are independent (permutation test), IHD dataset.} \label{Img:error_ind_ibs}
\end{figure}

Table \ref{tabl:ibs_best} contains (estimated) prediction errors of the most
significant combinations obtained by MDR analysis of ischemic heart
disease data.  At Figure \ref{Img:error_ind_ibs} a plot of empirical
distribution function of prediction error is given when the disease
is not linked with explanatory variables. We use here the simulated
samples $\overline{\xi}_b$ introduced in Section \ref{pred_alg},
with $b=1,\ldots,B$ where $B=100.$ One can see that out of these 100
simulations, the corresponding prediction error was not less than
$0.42.$ Note that Monte Carlo ${\sf p}$-value
\eqref{montecarlopvalue} of all three combinations is less than
$0.01$ (since their prediction errors are much less than $0.42$),
which is usually considered as a good performance.

\begin{table}[h]
\begin{center}
\begin{tabular}{|l|c|}
\hline
Factors  & Prediction error\\
\hline
 GpIa, FXIII, AH, HC & 0.231\\
 \hline
 Cx37, AH, HC & 0.238 \\
\hline
 GpIa,Cx37, AH, HC & 0.241\\
 \hline
 \end{tabular}
 \end{center}
 \caption{The most significant combinations obtained by MDR
analysis for IHD
data.}\label{tabl:ibs_best}
\end{table}

Table \ref{tabl:ibs_best_mdr_ind_rule} contains the results of MDRIR
method, which are similar to results of MDR method. However, it is worth mentioning
that MDRIR method allows to identify additional
combinations with prediction error around $0.24$.

\begin{table}[h]
\begin{center}
\begin{tabular}{|l|c|}
\hline
Factors & Prediction error\\
\hline
 FXIII,  FVII, AH, HC & 0.240\\
\hline
FXIII,  AH, HC & 0.242 \\
\hline
GpIa, Cx37, AH, HC & 0.247 \\
\hline
\end{tabular}
\end{center}
\caption{The most significant combinations  obtained by MDRIR
analysis of IHD data.} \label{tabl:ibs_best_mdr_ind_rule}
\end{table}

\vskip0,3cm
It follows from Tables \ref{tabl:ibs_best} and
\ref{tabl:ibs_best_mdr_ind_rule} that hypertension and
hypercholesterolemia are the most important external risk factors.
Indeed, these two factors appear in every of 6 combinations.

To perform a more precise analysis of influence of SNPs on IHD
provoking we analyze gene-gene interactions. We used two different
strategies. Namely, we applied MDR method to a subgroup of
individuals who are not subject to any of the  external risk factors
(i.e. to non-smokers without obesity and without
hypercholesterolemia, 51 cases and 97 controls). Another strategy is to
apply MDR method to the whole sample, but to take into account only
genetic factors rather than all factors. Table
\ref{tabl:ibs_best_genes} contains the most significant combinations
of SNPs and their prediction errors.
\vspace{3mm}

\begin{table}[h]
\begin{center}
\begin{tabular}{|l|c|c|}
\hline
Method & Genetic factors & Prediction\\&& Error\\
\hline
MDR on a subgroup of individuals &&\\who
do not have any risk factors & GpIa, Cx37 & 0.281\\
\hline
MDR method on the whole group&&\\ taking into account only genetic factors & GpIa, Cx37 & 0.343\\
\hline
\end{tabular}
\end{center}

\vspace{3mm}

\caption{Comparison of the most significant SNP combinations
obtained by two different ways of MDR analysis of IHD data.}
\label{tabl:ibs_best_genes}
\end{table}

\vspace{5mm}

It turned out that both methods  yield   similar results.
Combination of SNPs GpIa and Cx37 has the biggest influence on IHD.
Prediction error is about 0.28-0.34, and smaller error corresponds
to a {\it risk-free sample}. Moreover it follows from Tables
\ref{tabl:ibs_best_mdr_ind_rule} and \ref{tabl:ibs_best_genes}
 that prediction error significantly
dropped after additional exogenous factors were taken into account
(the error is 0.247 if additional external factors are taken into
account and 0.343 if not).

Thus based on ischemic heart disease data with the help of Tables
\ref{tabl:ibs_best}--\ref{tabl:ibs_best_genes} we can make the
following conclusions. Combination of two SNPs (GpIa and Cx37) and
two external factors (hypertension and hypercholesterolemia) has the
biggest influence on IHD. Also FXIII gives additional predictive
power if AH and HC are taken into account.

\vspace{7mm} {\bf Myocardial infarction} \vspace{3mm}

\begin{figure}[h]
\begin{center}
 \includegraphics[width=8cm]{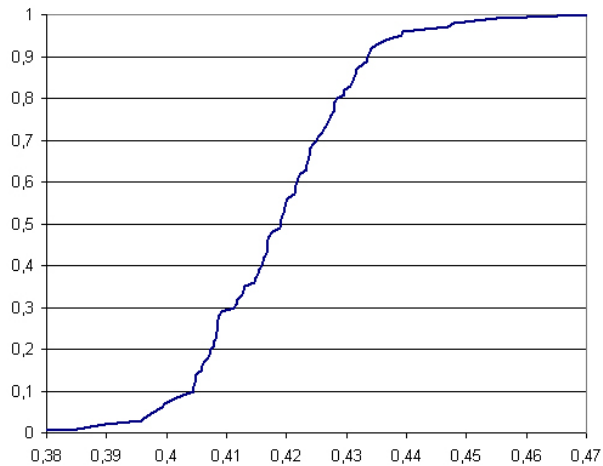}
   \end{center}
\caption{The empirical c.d.f.
of the estimated error for the case when disease risk and
predictors are independent (permutation test), MI dataset.} \label{Img:error_ind_im}
\end{figure}

Prediction errors of the most significant combinations obtained by
MDR analysis of MI data are presented in Table \ref{tabl:im_best}.
 Figure  \ref{Img:error_ind_im} contains the plot of
empirical
c.d.f.
of prediction error if disease is
not linked with risk factors. This curve shows that for all 100
simulations of $\overline{\xi}_b$ the estimated prediction error was
not less that 0.38. Note that Monte Carlo ${\sf p}$-value of all
combinations is less than
$0.01$.
 \vspace{3mm}
\begin{table}[h]
\begin{center}
\begin{tabular}{|l|c|}
\hline
Factors & Prediction error\\
\hline
 GpIIIa, FXIII, Cx37, AH & 0.343\\
 \hline
 GpIIIa, FXIII, FVII, Cx37 & 0.347 \\
\hline
 Cx37, Sm & 0.356\\
 \hline
 \end{tabular}
 \end{center}
 \caption{The most significant combinations obtained by MDR
analysis of MI dataset.}
 \label{tabl:im_best}
\end{table}


%

\vspace{3mm}
MDRIR analysis of the same dataset gives a  clearer picture (see
Table \ref{tabl:im_best_mdr_ind_rule}).

\begin{table}[!h]
\begin{center}
\begin{tabular}{|l|c|}
\hline
Factors & Prediction error\\
\hline
Cx37, Sm & 0.351\\
\hline
GpIIIa,  Cx37, Sm & 0.353 \\
\hline
GpIIIa,  Cx37, Sm, HC & 0.355 \\
\hline
\end{tabular}
\end{center}
\caption{The most significant combinations obtained by MDRIR
analysis of MI dataset.} \label{tabl:im_best_mdr_ind_rule}
\end{table}
Apparently, combination of smoking and SNP Cx37 is the most
significant. These two factors appear in all combinations in Table
\ref{tabl:im_best_mdr_ind_rule}. Involving any additional factors only
increases prediction error.

The explicit form of the prediction algorithm based on Cx37 and Sm
shows that these factors interact nonlinearly. Smoking as well as
Cx37 homozygote leads to the disease. However wild-type allele can
protect from consequences of smoking, because combination of smoking
and Cx37 wild-type is  a protective one (i.e. value of prediction
algorithm of this combination is -1).

\subsection{Logic regression}

We performed several research procedures both  for IHD and MI data,
with different restrictions imposed on the statistical model. To
describe these models set
$$(X_1,\dots,X_n)=(Z_1,\dots,Z_m,R_1,\dots,R_k)$$
where variables $Z=(Z_1,\dots,Z_m)$ stand for SNP values (PAI-1,
GpIa, GpIIIa, FXIII, FVII, IL-6, Cx37 respectively) and
$R=(R_1,\dots,R_k)$ denote external risk factors (Ob, AH, Sm, HC),
$m=7,$ $k=4.$ We consider four different models in order to analyze
both total influence of genetic and external factors and losses in
predictive force appearing when some factors are excluded.
In our
applications we will take $s=3,$ as search over larger forests for
samples with modest sizes can give very complicated and unreliable
results.
\vspace{2mm}

\textbf{Model 1.}  We consider the class $\mathcal{M}$ (see Section 2.3) consisting of
 the functions $h$ having a form
$$ h(Z,R)=\beta_0+\sum_{v=1}^s \beta_v
T_v(Z_1,\ldots,Z_{m})+\sum_{v=1}^k \beta_{s+v}R_v
$$ where the coefficients $\beta_i \in \mathbb{R}$ and $T_i$ are polynomials identified with trees.
In other words we require that external factors are present only in
trees consisting of one variable.

\textbf{Model 2.}   Now we assume that any function $h\in \mathcal{M}$
has the representation
\begin{equation}\label{model2} h(Z,R)=\beta_0+\sum_{v=1}^s \beta_v
T_v(Z_1,\ldots,Z_{m},R_1,\dots,R_k) \end{equation} where $\beta_i \in \mathbb{R}$ and $T_i$ are polynomials
identified with trees. Thus we allow the interaction of genes and external factors in order to find significant
gene-environment interactions.  However we impose additional restrictions to avoid too complex combinations of
external risk factors. We do not tackle here effects of interactions where several external factors are
involved. Namely, we consider only the trees satisfying the following two conditions.
\begin{enumerate}
\item If there is a leaf containing external factor variable then
the root of that leaf contains product operator.
\item Moreover, another branch growing from the same root is also a leaf and contains a genetic (SNP) variable.
\end{enumerate}


{\bf Models 3 and 4} have additional restrictions that polynomials
$T_v\,(v=1,\dots,s)$ in \eqref{model2} depend only on external
factors and only on SNPs respectively. These models are considered
to compare their results with ones obtained with all information
taken into account, in order to demonstrate the importance of
genetic (resp. external) data for risk analysis.

\vspace{5mm} {\bf Ischemic heart disease}
\vspace{2mm}

 We have the
following results.

\vspace{3mm}
\begin{table}[h]
\begin{center}
\begin{tabular}{|l|c|c|c|c|}
\hline
Model & 1 & 2 & 3 & 4 \\
\hline
Prediction error & 0.19040 & 0.20364 & 0.22812 & 0.33990\\
\hline
\end{tabular}
\end{center}\caption{Results of LR for IHD dataset.}\label{tabl:ibs_annealing}
\end{table}

Note that prediction error in Model 1 is only about 0.19. For the same model we performed also {\it fast
simulated annealing} search of the optimal forest which is much more time-efficient, and a reasonable error of
0.23 was obtained. Model 3 application shows that external factors play an important role in IHD genesis, as
classification based on external factors only gives the error less than  0.23, while usage of SNPs only   (Model
4) lets the error grow to 0.34.

Model 1 gave the minimal prediction error. For the optimal forest
$(T_1,\ldots,R_4)$ the function
  $\widehat{h}(Z,R)$ given before formula \eqref{lr_predalg2} with
  $S=\{1,\ldots,N\}$ is
provided by the expression
\begin{equation}\label{forh}
-0.597 T_1-0.354 T_2+ 0.521 T_3-0.444 R_1+1.311 R_2- 0.146
R_3 +2.331 R_4-0.226
\end{equation}
where\footnote{The sums and products are
modulo   $3$.}
$$
    T_1 = (Z_4 Z_3+ Z_6   Z_7+ Z_2   Z_2 + Z_3    Z_7)    (Z_1 )^2        Z_3
    Z_7,$$ $$
    T_2 = Z_1   (Z_3)^2         (Z_6   Z_7 + Z_7   (Z_4)^2
    Z_2),\;\;
    T_3 = Z_2 + 2Z_2  (Z_6)^2        Z_7.$$
\vspace{3mm}

The external factors  2  and 4 (i.e. AH and HC) are the most
influential since the coefficients at them are the greatest ones
(1.311 and 2.331). As was shown above,  MDR yields the same
conclusion. If the gene-environment interactions are allowed (Model
2), no considerable increase in predictive force has been detected.
However we list the pairs of SNPs and external factors present in
the best forest:  $Z_7$ and $R_2$, $Z_7$ and $R_1$, $Z_7$ and $R_4$,
$Z_5$ and $ R_1.$ It is seen that   Cx37 SNP is of substantial
importance as it appears in combination with all risk factors except
for smoking.

As formula \eqref{forh} is hard to interpret, we select the most
significant SNPs via a variant of permutation test. Consider a
random rearrangement of the column with first SNP in IHD dataset.
Calculate the prediction error using these new simulated data and
the same function $\widehat{h}$ as before. The analogous procedure
is done for other columns (containing the values of other SNPs) and
the errors found are given in Table \ref{tabl:ibs_vim}. It is seen
that the error increases considerably when the values of
 GpIa and Cx37 are permutated. The statement that they are the main sources
 of risk agrees with what was obtained above by MDR
 method.

\vspace{3mm}
\begin{table}[h]
\begin{center}
\begin{tabular}{|p{2.9cm}|p{1.3cm}|p{1.3cm}|p{1.3cm}|p{1.3cm}|p{1.3cm}|p{1.3cm}|p{1.3cm}|}
\hline
Prediciton error for Model 1& GpIa & Cx37 & IL-6 &  PAI-1  & GpIIIa & FXIII & FVII  \\
\hline
0.19040 & 0.26283 & 0.25987 & 0.22590 & 0.21212  & 0.20798 & 0.20173 & 0.19040  \\
\hline
\end{tabular}
\end{center}\caption{The SNP significance test for IHD in Model 1.}\label{tabl:ibs_vim}
\end{table}

\vspace{2mm} {\bf Myocardial infarction}
\vspace{2mm}

For the MI dataset, under the same notations that above, the
following results for our four models were obtained.

\vspace{3mm}
\begin{table}[h]
\begin{center}
\begin{tabular}{|l|c|c|c|c|}
\hline
Model   & 1 & 2 & 3 & 4 \\
\hline
Prediction error & 0.30526 & 0.33058 & 0.39057 & 0.36455\\
\hline
\end{tabular}
\end{center}\caption{Results of LR for IHD dataset.}\label{tabl:im_annealing}
\end{table}

To comment the Table \ref{tabl:im_annealing} we should first
underline that external risk factors play less important role
compared with IHD risk: if they are used without genetic
information, the error increases by 0.09, see Models 1 and 3 (while
the same increase for IHD was  0.03).
The function
  $\widehat{h}(Z,R)$ defined before formula \eqref{lr_predalg2} with
  $S=\{1,\ldots,N\}$ is
is equal to
$$-1.144
T_1+0.914 T_2- 0.45 T_3-0.285 R_1-0.675 R_2+ 0.828 R_3
-0.350R_4-0.055 $$  where $$
    T_1 = Z_1  Z_3  (Z_5)^2   ,\;\;
    T_2 = Z_7,\;\;
    T_3 = Z_4 + Z_3 + Z_7 + Z_6.
$$
Thus the first tree has the greatest weight  (coefficient equals
-1.144), the second tree (i.e.  Cx37 SNP) is on the second place,
and external factors are less important.

As for IHD we performed a permutation test to compare the
significance of different SNPs.
Its results are presented in Table \ref{tabl:im_vim}.

\vspace{3mm}
\begin{table}[h]
\begin{center}
\begin{tabular}{|p{2.9cm}|p{1.3cm}|p{1.3cm}|p{1.3cm}|p{1.3cm}|p{1.3cm}|p{1.3cm}|p{1.3cm}|}
\hline
Prediciton error for Model 1 & Cx37 & GpIIIa & IL-6 & FXIII & FVII  &  PAI-1 & GpIa    \\
\hline
0.30526 & 0.44420 & 0.35345 & 0.33998 & 0.32761 & 0.32427 & 0.31918 & 0.30526    \\
\hline
\end{tabular}

\end{center}\caption{The SNP significance test for MI in Model 1.}\label{tabl:im_vim}
\end{table}

As seen from this table, the elimination of  Cx37 SNP leads to a
noticeable increase in the prediction error. This fact agrees with
results obtained by MDR analysis of the same dataset.

\vspace{9mm}
\subsection{Results obtained by RF and SGB methods}
\vspace{2mm} The given datasets were unbalanced w.r.t. response
variable and we first applied the resampling technique to them. That
is, enlargement of the smaller of two groups case-control in the
sample by additional bootstrap observations till the final
proportion case:control would be  1:1. We also employed
modifications of RF by \cite{Breiman04} and SGB by \cite{Liu08} for
unbalanced samples, but those worked poorly for permutation tests
and we do not give their results here. Note that due to the
resampling techniques the following effect arise. Some observations
in small groups (case or control) appear in the new sample more
frequently than other ones. Therefore, we took the average over 1000
iterations.

\vspace{5mm} {\bf Ischemic heart disease}
\vspace{2mm}

\begin{table}[h]
\begin{center}{
\begin{tabular}{|l|c|c|c|c|c|c|c|c|c|c|c|c|}\hline
Data &RF&SGB\\
\hline
with SNP  &0.20/0.454&0.134/0.473\\
\hline
without SNP  &0.23/0.51&0.261/0.503\\
\hline
\end{tabular}}
\end{center}
\caption{ \small Prediction error$/$prediction error in permutation
test calculated via cross-validation for IHD dataset with employment
of RF  and SGB methods.\label{tabl:sgb1}}
\end{table}

Results of RF and SGB methods are given in  Table
\ref{tabl:sgb1}. It shows that RF and SGB methods give statistically
reliable results (prediction error in the permutation test is
close to $50\%$). Moreover, additional SNP information improves
predicting ability on 11$\%$ and 13$\%$ (SGB). It seems that SGB method is better
fitted to IHD data than RF.

Computing CVIM for each $X_i$,  we constructed $Z_i$ as follows. We
included in $Z_i$ all predictors $X_j,$ $j\neq i$, for which
$\chi^2$-criteria rejected hypothesis of independence between $X_j$
and $X_i$ at $5\%$ significance level. Since the genetic information
has second order effect on prediction of $Y$ comparing to the risk
factors, we ran the program 1000 times and then took the average
CVIM to get a reliable estimate. An error over different runs of the
program was around 0.01. The results are given in
Table 11.
\vspace{2mm}
\begin{table}[h]
\begin{center}{
\begin{tabular}{|c|c|c|c|c|c|c|c|c|c|c|c|c|}\hline
AH &HC&Cx37 &Ob &FXIII &Sm &GpIa &FVII&PAI-1& GpIIIa  &IL-6  \\
\hline
8.9&5.3&5.1&0.56& 0.53&0.11& 0.1& 0.07& 0.03&0.02&0.01\\
\hline
\end{tabular}}\end{center}
\caption{\small
  Predictors are ranged in terms of their CVIM  for IHD dataset.}
\end{table}

Thus, the most relevant predictors for IHD are AH, HC and Cx37.

{\vspace{5mm}} {\bf Myocardial infarction}
\vspace{2mm}

Results of RF and SGB methods are given in the following table.

\begin{table}[h]
\begin{center}
\begin{tabular}{|l|c|c|c|c|c|c|c|c|c|c|c|c|}\hline
Data&RF&SGB\\
\hline
with SNP data&0.36/0.497&0.399/0.53\\
\hline
without SNP data&0.473/0.527&0.482/0.562\\
\hline
\end{tabular}
\end{center}
\caption{ \small Prediction error$/$prediction error in permutation
test calculated via cross-validation for MI dataset with employment
of RF  and SGB methods.\label{tabl:sgb2}}
\end{table}
Table \ref{tabl:sgb2} shows that RF and SGB methods give
statistically reliable estimates (prediction error in the
permutation test is close to $50\%$). Moreover, additional SNP
information improves predicting ability on 10$\%$.

CVIM was calculated according to \eqref{cvim} and is given below.

\begin{table}[h]
\begin{center}{
\begin{tabular}{|c|c|c|c|c|c|c|c|c|c|c|c|c|}\hline
 Cx37& Sm &AH & GpIIIa &FVII& FXIII& HC& GpIa& Ob& IL-6&PAI-1\\
 \hline
 7.5&2&1.86&0.03&0.02& $\approx$0&$\approx$0&$\approx$0&$\approx$0& $\approx$0&$\approx$0\\
\hline
\end{tabular}}
\end{center}
\caption{ \small  Predictors are ranged in terms of their CVIM for  MI
dataset.}
\end{table}
Thus, the most relevant predictors for MI are  Cx37, Sm and AH.

\section{Conclusions and final remarks}

Let us briefly summarize  the main results obtained. The analysis of
IHD dataset showed that two external risk factors out of four
considered (AH and HC) have a strong connection with the disease
risk  (the error of classification based on external factors only is
0.25--0.26 with ${\sf p}$-value less than  0.01). Also, the classification
based on SNPs only gives  a relatively low error of  0.28. Moreover,
the most influential SNPs are   Cx37 and GpIa (FXIII also enters the
analysis only when AH and HC are present). Prediction error
decreases to
  0.13 if  both
 SNP information and external risk factors
 are taken into account. Note that excluding any of the 5 remaining SNPs (all except for two most influential) from data increases the error by
  0.01--0.02 approximately. So, while the most influential data are
  responsible for the situation within  a large part of population,
  there are smaller parts where other SNPs come to effect and
  provide a more efficient prognosis  (``small subgroups effect$"$).

 The MI dataset gave the following results. The most significant factors of MI risk are the
 Cx37 SNP (more precisely, homozygous mutation) and smoking with a considerable   gene-environment interaction
 present.
The smallest prediction error of methods applied was 0.33--0.35
(with ${\sf p}$-value less than  0.01). The classification based on
external factors only yields a much greater error of 0.42. Thus
genetic data improves the prognosis quality essentially. While two
factors are of great importance, other SNPs considered actually do
not improve the prognosis essentially, i.e. no small groups effect
is observed.

The conclusions given above are based on several complementary
methods of modern statistical analysis. These new data mining
methods allow to analyze other datasets as well. The study can be
continued with larger datasets, in particular,  involving new SNP
data.

\newpage
{\Large {\bf Acknowledgments}}
\vskip0,4cm

$\!\!$The authors are grateful to Professor V.A.Tkachuk, Associate Professor
L.M.Samokhod\-skaya and MD. A.V.Balatsky for providing the data concerning the complex
diseases. On the basis of this preprint the joint paper
will be prepared with medical interpretation
of the obtained results.

$\!\!$A.V.Bulinski also would like to thank Professors I.Kourkova and G.Pag\`es for
invitation to LPMA of the University Pierre and Marie Curie, he is grateful
to all the members of LPMA for hospitality.

\vspace{3mm}
$\!\!$The work is partially supported by RFBR grant 10-01-00397.

\newpage

{\bf Alexander BULINSKI},
\vskip0,4cm
Faculty of Mathematics and Mechanics, Lomonosov Moscow State University,

GSP-1, Leninskie gory, Moscow, 119991, Russia
\vskip0,2cm
{\it and}
\vskip0,2cm
LPMA UPMC University Paris-6,

4 Place Jussieu, 75252 Paris CEDEX 05, France

\vskip0,4cm
{\it E-mail address}: bulinski@mech.math.msu.su

\vskip1,2cm
{\bf Oleg BUTKOVSKY},
\vskip0,4cm
Faculty of Mathematics and Mechanics, Lomonosov Moscow State University,

GSP-1, Leninskie gory, Moscow, 119991, Russia
\vskip0,4cm
{\it E-mail address}: oleg.butkovskiy@gmail.com

\vskip1,2cm
{\bf Alexey SHASHKIN},
\vskip0,4cm
Faculty of Mathematics and Mechanics, Lomonosov Moscow State University,

GSP-1, Leninskie gory, Moscow, 119991, Russia
\vskip0,4cm
{\it E-mail address}: ashashkin@hotmail.com

\vskip1,2cm
{\bf Pavel YASKOV},
\vskip0,4cm
Steklov Mathematical Institute, Gubkina str. 8, Moscow, 119991, Russia
\vskip0,2cm
{\it and}
\vskip0,2cm
Faculty of Mathematics and Mechanics, Lomonosov Moscow State University,

GSP-1, Leninskie gory, Moscow, 119991, Russia
\vskip0,4cm
{\it E-mail address}: pavel.yaskov@mi.ras.ru


\newpage

\begin{thebibliography}{99}

\bibitem{Albrechtsen}
A. Albrechtsen, S. Castella, G. Andersen, T. Hansen, O. Pedersen and
R. Nielsen. {\it A Bayesian multilocus association method: allowing
for higher-order interaction in association studies.} Genetics, vol.
176 (2007),  pp. 1197-1208.

\bibitem{ArlotCelisse} A. Arlot and A. Celisse. {\it A survey
of cross-validation procedures for model selection.} Statist. Surv.,
vol.4 (2010), pp. 40-79.

\bibitem{Bahadur61}
R. Bahadur. {\it A representation of the joint distribution of
responses to $n$ dichotomous items.} Studies in Item Analysis and
Prediction, Stanford University Press, H. Solomon (ed.), 1961, pp.
158-168.

\bibitem{Biau10}
G. Biau. {\it Analysis of a Random Forests model.} J. of Machine
Learning Research, LSTA, LPMA, Universit\'e Paris-6, 2010.

\bibitem{Biau08}
G. Biau, L. Devroye and G. Lugosi. {\it Consistency of Random
Forests and Other Averaging Classifiers.} J. of Machine Learning
Research, vol.9 (2008), pp. 2015-2033.

\bibitem{Bock} C. Bock   and T. Lengauer. {\it Computational
epigenetics.} Bioinformatics, vol. 24 (2008), pp. 1-10.

\bibitem{Bondy} A. Bondy and U.S.R. Murty. {\it Graph Theory.}
Springer, 2008.

\bibitem{Breiman04} L. Breiman, C. Chen and A. Liaw. {\it Using random forest to learn imbalanced
data.} J. of Machine Learning Research, no. 666, Department of
Statistics, University of California, Berkeley, CA, 2004.

\bibitem{Bulinskiconference} A. Bulinski. {\it Stochastic methods of identification of SNP interactions.}
The 1st Int. Research and Practice Conference on Postgenomic Methods
of Analysis in Biology, and Laboratory and Clinical Medicine, MSU,
2010, p. 146.

\bibitem{ourreport}
A. Bulinski, O. Butkovsky, A. Shashkin, P. Yaskov, M. Atroshchenko
and A. Khaplanov. {\it Investigations in the framework of the MSU
Research project Postgenomic Medical Studies and Technologies.} MSU,
2010.

\bibitem{Chawla10} N.V. Chawla. {\it Data mining for imbalanced datasets: An overview.}
Data mining and knowledge discovery handbook, Part 6, O. Maimon and
L. Rokach (eds.), Springer, 2010, pp. 875-886.

\bibitem{Coffey} C.S. Coffey, P.R. Hebert, M.D. Ritchie, H.M. Krumholz et
al. {\it An application of conditional logistic regression and
multifactor dimensionality reduction for detecting gene-gene
interactions on risk of myocardial infarction: the importance of
model validation.} Atherosclerosis, vol. 30 (2004), pp. 5-49.

\bibitem{Logit} D.R. Cox. {\it The analysis of multivariate binary data.}
Applied Statistics, vol. 21 (1972), pp. 113-120.

\bibitem{Dai} J.Y.Dai et al. {\it SHARE: an adaptive algorithm to select the most informative set of
SNPs for candidate genetic association.} Biostatistics, vol. 10
(2009), pp. 680-693.

\bibitem{Friedman01} J.H. Freedman. {\it Greedy function approximation: A gradient boosting
machine.} Ann. Statist., vol. 29 (2001), pp. 1189-1232.

\bibitem{FriIckst} A. Fritsch  and K. Ickstadt. {\it Comparing logic regression based
methods for identifying SNP interactions.} Lecture Notes in Computer
Science 4414, Springer, 2007, pp. 90-103.

\bibitem{GAW} GAW16,
http://www.ncbi.nlm.nih.gov\-/projects/gap/cgi-bin/study.cgi\-$?$study$\_$\-id=phs000128.v3.p3.

\bibitem{Golland05} P. Golland, F. Liang, S. Mukherjee and D.
Panchenko. {\it Permutation tests for classification}, Lecture notes
in Computer Science, Springer, vol. 3559 (2005), pp. 501-515.

\bibitem{Hajek} B. Hajek. {\it Cooling schedules for optimal annealing.}
Math. Oper. Res., vol. 13 (1988), pp. 311-329.

\bibitem{HapMap} HapMap, http://www.hapmap.org.

\bibitem{Hastie10} T. Hastie, R. Tibshirani and J. Friedman. {\it The Elements of
Statistical Learning: Data Mining, Inference, and Prediction.}
Springer, 2009.

\bibitem{logreg} D.  Hosmer and S. Lemeshow. {\it Applied Logistic
Regression.} Wiley, 2000.

\bibitem{Karchin} R. Karchin. {\it Next generation tools for the annotation of human SNPs.}
Briefings in Bioinformatics, vol. 10 (2009), pp. 35-52.

\bibitem{Tabara} K. Kohara, Y. Tabara, J. Nakura, Y.Imai et al. {\it Identification of
Hypertension-Susceptibility Genes and Pathways by a Systemic
Multiple Candidate Gene Approach: The Millennium Genome Project for
Hypertension.} Hypertension Research, vol. 31 (2008), pp. 203-212.

\bibitem{KooperBis} C. Kooperberg, J.C. Bis, K.D. Marciante, S.R. Heckbert, T. Lumley and B.M.
Psaty. {\it Logic regression for analysis of the association between
genetic variation in the renin-angiotensin system and myocardial
infarction or stroke.} Am. J. of Epidemiology, vol. 165 (2007), pp.
334-343.

\bibitem{Lee} S.Y. Lee, Y. Chung, R.C. Elston, Y. Kim and T. Park.
{\it Log-linear model-based multifactor dimensionality reduction
method to detect gene-gene interactions.} Bioinformatics, vol. 23
(2007), pp. 2589-2595.

\bibitem{Lehmann} E.L. Lehmann and J.P. Romano. {\it Testing Statistical
Hypotheses}, Springer, 2005.

\bibitem{Lengauer07} T. Lengauer (ed.). {\it Bioinformatics – From Genomes to
Therapies}, Wiley – VCH Verlag GmbH and KGaA, Weihhein, 2007.

\bibitem{Liu08} X.-Y. Liu,  J. Wu and Z.-H. Zhou. {\it Exploratory Undersampling for Class-Imbalance
Learning.} Systems, Man, and Cybernetics, Part B: Cybernetics, IEEE
Transactions, vol. 39 (2008), pp. 539-550.

\bibitem{Massart03} P. Massart. {\it Concentration Inequalities and Model
Selection.} Springer, 2003.

\bibitem{MotRit2006} A.A. Motsinger and M.D. Ritchie. {\it The effect of reduction in
cross-validation intervals on the performance of multifactor
dimensionality reduction.} Genetic Epidemiology, vol. 30 (2006), pp.
546-555.

\bibitem{Nikolaev} A.G. Nikolaev and  S.H. Jacobson. {\it Simulated
Annealing.} Handbook of Metaheuristics, Springer, 2010, pp. 1-39.

\bibitem{Nunke} R. Nunkesser, T. Bernholt, H. Schwender, K. Ickstadt and I.
Wegener. {\it Detecting high-order interactions of single nucleotide
polymorphisms using genetic programming.} Bioinformatics, vol. 23
(2007), pp. 3280-3288.

\bibitem{Olive} D.J. Olive, {\it  The number of samples for resampling
algorithms.} http://www.math.siu.edu/olive/ppresamp.pdf, 2010.

\bibitem{Park09} J. Park. Independent rule in classification of multivariate binary
data. J. Multivar. Anal., vol. 100 (2009), pp. 2270-2286.

\bibitem{MDR} M.D. Ritchie, L.W. Hahn, N. Roodi, R. Bailey, W.D. Dupont, F.F. Parl and
J.H.Moore. {\it Multifactor-dimensionality reduction reveals
high-order iteractions among estrogen-metabolism genes in sporadic
breast cancer.} Am. J. Hum. Genet., vol. 69 (2001), pp. 138-147.

\bibitem{Sylvain09} S. Rubenthaler, T.  Ryden  and M.   Wiktorsson.
{\it Fast simulated annealing in $\mathbb{R}^d$ with an application
to maximum likelihood estimation in state-space models.} Stoch.
Proc. Appl., vol. 119 (2009), pp. 1912-1931.

\bibitem{RuczKB} I. Ruczinski, C. Kooperberg and M.  LeBlanc.
{\ it Logic regression.} J. Comp.  Graph. Statist., vol. 12 (2003),
pp. 475-511.

\bibitem{Schwender}
H. Schwender and K. Ickstadt. {\it Identification of SNP
interactions using logic regression.} Biostatistics, vol. 9 (2008),
pp. 187-198.

\bibitem{Schwender2} H. Schwender and K. Ickstadt.
{\it Empirical Bayes analysis of single nucleotide polymorphisms.}
Bioinformatics, vol. 9 (2008), pp. 144-159.

\bibitem{Schwender3}
H. Schwender and I. Ruczinski. {\it Testing SNPs and sets of SNPs
for importance in association studies.} Biostatistics, vol. 11
(2010), pp. 1-15.

\bibitem{Strobl08} C. Strobl, A.L. Boulesteix, T. Kneib, T. Augustin and
A. Zeileis. {\it Conditional variable importance for random
forests.} BMC Bioinformatics, vol. 9 (2008), p. 307.

\bibitem{MachLearn09} S. Szymczak, J.M. Biernacka, H.J. Cordell, G.-R. Oscar, I.R. K\"onig, H. Zhang
and Y.V. Sun. {\it Machine Learning in Genome-Wide Association
Studies.} Genetic Epidemiology, vol. 33 (2009), pp. 51-57.

\bibitem{ChiaTiTsai} C.-T. Tsai, J.-J. Hwang, M.D. Ritchie and J.H. Moore.
{\it Renin–angiotensin system gene polymorphisms and coronary artery
disease in a large angiographic cohort: Detection of high order
gene$-$gene interaction.} Atherosclerosis, vol. 195 (2007), pp.
172-180.

\bibitem{Velez} D. Velez, B.C. White, A.A. Motsinger, W.S. Bush, M.D. Ritchie, S.M. Williams and J.H.
Moore. {\it A balanced accuracy function for epistasis modeling in
imbalanced datasets using multifactor dimensionality reduction.}
Genetic Epidemiology, vol. 31 (2007), pp. 306-315.

\bibitem{SNPgeneral} J. Venter, M.D. Adams, E.W. Myers, P.W. Li et al.
{\it The Sequence of the Human Genome.} Science, vol. 291 (2001),
pp. 1304-1351.

\bibitem{SNPHunter09} X. Wan, C. Yang, Q. Yang, H. Xue, N.L.S. Tang and W.
Yu. {\it MegaSNPHunter: a learning approach to detect disease
predisposition SNPs and high level interactions in genome wide
association study.} BMC Bioinformatics, vol. 10 (2009), p. 13.

\bibitem{CrossValid} S. Winham, A.J. Slater and A.A. Motsinger-Reif. {\it A comparison of
internal validation techniques for multifactor dimensionality
reduction.} BMC Bioinformatics, vol. 11 (2010), pp.394.

\bibitem{ZhangBuc} Z. Zhang, E.S. Buckler, T.M. Casstevents and P.J.
Bradburg. {\it Software engineering the mixed model for genome wide
assosiation studies on large samples}, Briefings in Bioinformatics,
vol. 10 (2009), pp. 664-675.








\end{thebibliography}
\end{document}